\documentstyle[floats,prl,aps]{revtex}
\begin{document}
\renewcommand{\thefootnote}{\fnsymbol{footnote}}
\draft
                     


\def\a{\alpha}
\def\b{\beta}
\def\d{\delta}
\def\e{\epsilon}
\def\g{\gamma}
\def\k{\kappa}
\def\l{\lambda}
\def\o{\omega}
\def\t{\theta}
\def\s{\sigma}
\def\D{\Delta}
\def\L{\Lambda}
\def\A{U_q(A^{(2)}_2)} 
\def\Uq{U_q(\widehat{sl}(2|1)}
\def\R{\overline{R}}

\def\beq{\begin{equation}}
\def\eeq{\end{equation}}
\def\bea{\begin{eqnarray}}
\def\eea{\end{eqnarray}}
\def\ba{\begin{array}}
\def\ea{\end{array}}
\def\no{\nonumber}
\def\le{\langle}
\def\re{\rangle}
\def\lt{\left}
\def\rt{\right}

\newtheorem{Theorem}{Theorem}
\newtheorem{Definition}{Definition}
\newtheorem{Proposition}{Proposition}
\newtheorem{Lemma}{Lemma}
\newtheorem{Corollary}{Corollary}
\newcommand{\proof}[1]{{\bf Proof. }
        #1\begin{flushright}$\Box$\end{flushright}}

\newcommand{\sect}[1]{\setcounter{equation}{0}\section{#1}}
\renewcommand{\theequation}{\thesection.\arabic{equation}}

\title{\large\bf The twisted quantum affine algebra $\A$ and
correlation functions of the Izergin-Korepin model}
\author{\large Bo-Yu Hou$^{~1}$, Wen-Li Yang$^{~1,~2}$
and Yao-Zhong Zhang$^{~2}$}

\address{$^{1}~$ Institute of Modern Physics,
Northwest University, Xian  710069,China\\
$^{2}~$ Department of Mathematics, University of Queensland,
Brisbane, Qld 4072, Australia}

\maketitle
\vspace{10pt}

\begin{abstract}
We derive the exchange relations of the
vertex operators of $U_q(A_2^{(2)})$ and show that these vertex operators
give the bosonization of the Izergin-Korepin model. We give an integral
expression of the correlation functions of the Izergin-Korepin model 
and derive the difference equations which they satisfy.
\end{abstract}

\section{Introduction} 
It is well-known that quantum affine algberas play an essential
role in the studies of low-dimensional massive integrable
models such as  quantum spin chains, since they provide the symmetry
algberas of these models.Based on the works of
 q-vertex operators\cite{Fre92}\cite{Fre88}
and level-one
highest weight representations of the quantum affine algebra
$U_q(\widehat{sl}_2)$, the Kyoto group\cite{DFJMN}\cite{JM}
developed a new method
 which enabled the group to diagonalize the XXZ spin-$\frac{1}{2}$ chain
 directly in the thermodynamic limit and moreover compute the
 correlation functions of the spin operators.This approach was later
 generalized to higher spin XXZ chains \cite{ID,BW,HK}, vertex models with
 $U_q(\widehat{sl}_n)$ symmetries\cite{K}, the face type
statistical mechanics models \cite{LP,AJ}, and more recently to
 integrable models with quantum affine superalgbera symmetry \cite{YZ}.

In this paper, we will extend the above programme further to the
Izergin-Korepin  19-vertex model\cite{IK}---\cite{Mar} 
which has twisted quantum affine algebra $\A$ as
its non-abelian symmetry.The bosonic
realiztion of $\A$ and its level-one vertex operators were
constructed recently \cite{Jing}. In section 3, 
we use the results of \cite{Jing} to calculate  the
exchange relations of the $q$-vertex operators .We show that these
vertex operators satisfy the Faddeev-Zamolodchikov algebra with the
R-matrix  of $\A$ as its constructure constant.
A Miki's construction of $\A$ is also given. In section 4, 
generalizing the Kyoto group's work
\cite{DFJMN} to the case of twisted quantum affine algebra $\A$, 
we give the bosonization of the
Izergin-Korepin model .
In section 5,
 we compute the 
correlation functions of the local operators (such as the spin operator
$S_z$ ) and give an integral
expression of the correlation functions. A set of 
difference equations satisfied by the correlation functions have also
been derived.

\section{Preliminaries}

In this section, we briefly review the bosonization of the twisted quantum
affine  algebra $U_q(A^{(2)}_2)$ at level one \cite{Jing}.

\subsection{ Quantum affine algebra $U_q(A^{(2)}_2)$ }

The symmetric Cartan matrix of the twisted affine Lie algebra $A^{(2)}_2$
is 
\begin{eqnarray*}
(a_{ij})=\left(
\begin{array}{cc}
8&-4\\
-4&2
\end{array}\right)
\end{eqnarray*}
where $i,j=0,1$. Quantum affine algebra 
$U_q(A^{(2)}_2)$ is a $q$-analogue of the universal 
enveloping algebra of $A^{(2)}_2$ generated by the Chevalley 
generators $\{e_i,f_i,t_i^{\pm 1},d | i=0,1\}$, where $d$ is the
usual derivation operator.
The defining relations are\cite{CP}\cite{DGZ}
\begin{eqnarray*}
& & t_it_j =t_jt_i,\ \  t_id=dt_i, \ \ [d,e_i]=\delta_{i,0}e_i,\ \ 
[d,f_i]=-\delta_{i,0}f_i,\\ 
& &t_ie_jt_i^{-1}=(q^{1\over 2})^{a_{ij}}e_j,
\ \ t_if_jt_i^{-1}=(q^{1\over 2})^{-a_{ij}}f_j ,\\
& &[e_i,f_j] =\delta_{ij} \frac{t_i-t_i^{-1}}{q_i-q_i^{-1}},~~\\
& & \sum_{r=0}^{1-a_{ij}}(-1)^r
\left[
\begin{array}{c}
1-a_{ij}\\
r
\end{array}
\right]_{q_i}
(e_i)^re_j(e_i)^{1-a_{ij}-r}=0~~~,~~~{\rm if}~i\neq j,\\
& & \sum_{r=0}^{1-a_{ij}}(-1)^r
\left[
\begin{array}{c}
1-a_{ij}\\
r
\end{array}
\right]_{q_i}
(f_i)^rf_j(f_i)^{1-a_{ij}-r}=0~~~,~~~{\rm if}~i\neq j,\\
\end{eqnarray*}
\noindent where $~q_1=q^{1\over 2}~, ~q_0=q^2$ , $t_i=q_i^{h_i}$ , and 
\begin{eqnarray*}
& &[n]_{q_i}=\frac{q_i^n-q_i^{-n}}{q_i-q_i^{-1}}~~~,
~~~[n]_{q_i}!=[n]_{q_i}[n-1]_{q_i}\cdot\cdot\cdot [1]_{q_i}~~,\\
& &\left[
\begin{array}{c}
n\\
r
\end{array}
\right]_{q_i}=\frac{[n]_{q_i}!}{[n-r]_{q_i}!
[r]_{q_i}!}
\end{eqnarray*}

$U_q(A^{(2)}_2)$ is a quasi-triangular Hopf algebra
endowed with  Hopf algebra structure:
\bea
&&\Delta(t_i)=t_i\otimes t_i,~~~~
   \Delta(e_i)=e_i\otimes 1+t_i\otimes e_i,~~~~
    \Delta(f_i)=f_i\otimes t_i^{-1}+1\otimes f_i ,\\
\label{counit}
&&\epsilon(t_i)=1,~~~~\epsilon(e_i)=\epsilon(f_i)=0,\\
&&S(e_i)=-t_i^{-1} e_i,~~~~ S(f_i)=-f_i t_i, ~~~~
   S(t_i^{\pm 1})=t_i^{\mp 1},~~~~ S(d)=-d.
\eea 

$\A$ can also be realized by 
the Drinfeld generators \cite{Jing}\cite{Dri88} $\{d,\ \ X^{\pm}_{m}$,
$a_n$, $K^{\pm 1},\gamma^{\pm 1/2} |  m \in 
{\bf Z}, n \in {\bf Z}_{\ne 0}\}$. 
The relations read
\begin{eqnarray}
\label{DRB1}
& &\gamma\ \  {\rm is\ \  central },\ \
[K,a_n]=0,\ \ [d,K]=0,\ \ [d,a_n]=na_n,\\
& &[a_m,a_n] =\delta_{m+n,0} 
{\{[4n]_{q_1}-(-1)^n[2n]_{q_1}\}(\gamma^m-\gamma^{-m})
\over m(q_1-q_1^{-1})} ,\\
& &KX^{\pm}_m =q^{\pm 1}X^{\pm}_m K,\ \
[d,X^{\pm}_m]=mX^{\pm}_m ,\\ 
& &[a_m,X^{\pm}_n]=\pm {\{[4m]_{q_1}-(-1)^m[2m]_{q_1}\} \over m}
\gamma^{\mp |m|/2}X^{\pm}_{n+m},\\ 
& &[X^{+}_m,X^{-}_n]=\frac{1}{ q_1-q_1^{-1}}
(\gamma^{(m-n)/2}\psi^{+}_{m+n} -\gamma^{-(m-n)/2 }
\psi^{-}_{m+n}),\\
\label{DRB2}
& &(z-wq^{\pm 2})(z+wq^{\mp 1})X^{\pm}(z)X^{\pm}(w)=
(zq^{\pm 2}-w)(zq^{\mp 1}+w)X^{\pm}(w)X^{\pm}(z)~~,
\end{eqnarray}
\noindent where the corresponding Drinfeld currents
$\psi^{\pm}(z)$ and $X^{\pm}(z)$ are defined by
\begin{eqnarray*}
& &\psi^+(z)=\sum_{m=0}^{\infty}\psi^{+}_mz^{-m}=
Kexp\{(q_1-q_1^{-1})\sum_{k=1}^{\infty}a_kz^{-k}\}~~,\\
& &\psi^-(z)=\sum_{m=0}^{\infty}\psi^{-}_{-m}z^{m}=
K^{-1}exp\{-(q_1-q_1^{-1})\sum_{k=1}^{\infty}a_{-k}z^{k}\}~~~,\\
& &X^{\pm}(z)=\sum_{n\in Z}X^{\pm}_mz^{-m}~~.
\end{eqnarray*}

The Chevalley generators are related to the Drinfeld generators
by the formulae:
\begin{eqnarray}
\label{CB1}
& &t_1= K,\ \ e_1 = X^{+}_0,\ \ t_0=\gamma K^{-2},
\ \ f_1 = X^{-}_0 ,\\  
\label{CB2}
& &e_0= K^{-2}[X_0^{-},X^{-}_1]_q,~~~
f_0=\frac{1}{[4]_{q_1}^2}[X^{+}_{-1},X^{+}_0]_{q^{-1}}K^2~~.
\end{eqnarray}

\subsection{Bosonization of $\A$ at level one}
Let us introduce the bosonic $q$-oscillators\cite{Jing}
$\{a_n~,Q~, P| n\in Z-\{0\}\}$  which satisfy the commutation relations
\begin{eqnarray*}
& &[a_m,a_n]=\delta_{m+n,0}\frac{1}{m}
\{[4m]_{q_1}-(-1)^m[2m]_{q_1}\}
[2m]_{q_1},\\
& &[P,a_m]=[Q,a_m]=0~~,~~[P,Q] = 1,
\end{eqnarray*}
Set $ F_1=\oplus_{n\in Z}C[a_{-1},a_{-2},.....]e^{\frac{Q}{2}+nQ}|0>$,
where the Fock vacuum vector $|0>$ is defined by
\begin{eqnarray*}
a_n|0>=0~~~{\rm for}~ n>0~~~,~~~~P|0>=0.
\end{eqnarray*}
\noindent Then, on $F_1$ the
Drinfeld currents of $\A$ at level-one are realized 
by the free boson fields as \cite{Jing}
\begin{eqnarray}
& & \gamma=q~~,~~ K=q^P~~~,\no\\
& &\psi^{+}(z)=q^Pexp\{(q_1-q_1^{-1})\sum_{n=1}^{\infty}a_nz^{-n}\}~~,\no\\
& &\psi^{-}(z)=q^{-P}exp\{-(q_1-q_1^{-1})\sum_{n=1}^{\infty}a_{-n}z^{n}\}
~~~,\no\\
\label{CUR}
& &X^{\pm}(z) =
exp\{\pm\sum_{n=1}^{\infty}\frac{a_{-n}}{[2n]_{q_1}}q^{\mp\frac{n}{2}}z^n\}
exp\{\mp\sum_{n=1}^{\infty}\frac{a_{n}}{[2n]_{q_1}}q^{\mp\frac{n}{2}}z^{-n}\}
e^{\pm Q}z^{\pm P+\frac{1}{2}}~~.
\end{eqnarray}

By the relations between Chevalley basis and Drinfeld basis
(\ref{CB1}---\ref{CB2}), one can check that
\begin{eqnarray*}
& &t_1e^{\frac{Q}{2}}|0>=q^{\frac{1}{2}}e^{\frac{Q}{2}}|0>
=(q^{\frac{1}{2}})^{<\alpha_1,\Lambda_1>}e^{\frac{Q}{2}}|0>,\\
& &t_0e^{\frac{Q}{2}}|0>=e^{\frac{Q}{2}}|0>
=(q^{\frac{1}{2}})^{<\alpha_0,\Lambda_1>}e^{\frac{Q}{2}}|0>,\\
& &e_ie^{\frac{Q}{2}}|0>=0~~~, ~~i=0,1,
\end{eqnarray*}
\noindent where $\alpha_1,~\alpha_0=-2\alpha_1+\delta$ are the simple
roots of $A^{(2)}_2$, $\Lambda_1$ is one basical fundamental weight, and
$\delta$ is the imaginary root. The Fock space $F_1$ coincides with the
level-one irreducible 
highest weight module $V(\Lambda_1)$ with the highest weight vector
 given by $|\Lambda_1>=e^{\frac{Q}{2}}|0>$.

\subsection{Level-one vertex operators}
Let $V$ be  the 3-dimensional  evaluation  representation of
$\A$ ,$\{v_1, v_0, v_{-1}\}$  be  the basis vectors of $V$, and
 $E_{i,j}$ be the $3\times 3$ matrix whose $(i,j)$-element is unity
and zero  otherwise.Then
the 3-dimenional level-0 representation $V_z$ of $\A$ is given by 
\begin{eqnarray*}
& &e_1=\alpha^{-1}[2]_{q_1}E_{1,0}+\alpha E_{0,-1},~~~
e_0=zE_{-1,1},\\
& &f_1=\alpha^{-1}[2]_{q_1}E_{-1,0}+\alpha E_{0,1}~~,~~
f_0=z^{-1}E_{1,-1},\\
& & t_1=q_1^2E_{1,1}+q_1^{-2}E_{-1,-1}+E_{0,0},~~~~
t_0=q_0E_{-1,-1}+q_0^{-1}E_{1,1}+E_{0,0}~~,
\end{eqnarray*} 
\noindent where $\alpha=\{[2]_{q_1}
q^{-\frac{1}{2}}\}^{\frac{1}{2}}$.

We  define the dual modules $V_z^{*b}$ of $V_z$ by
$\pi_{V^{*b}}(a)=\pi_{V}(b(a))^{t}$, $\forall a\in \A$, where 
$t$ is the transposition operation. $b$ is the
anti-automorphism  defined by
\begin{eqnarray}
\label{DEF}
b(x)=(-q)^{\hat{\rho}}S(x)(-q)^{-\hat{\rho}}\stackrel{def}{=}
(-q)^{3d}(-q^{\frac{1}{2}})^{\frac{1}{2}h_1}S(x)
(-q)^{-3d}(-q^{\frac{1}{2}})^{-\frac{1}{2}h_1},~~~ \forall x\in \A ~,
\end{eqnarray}
\noindent where $2\hat{\rho}=6d+\frac{1}{2}h_1$.
A convenient
feature of $b$ is that $b^2=id$ (since 
 $S^2(x)=q^{-2\hat{\rho}}xq^{2\hat{\rho}}$) and that
the following isomorphism holds:
\begin{eqnarray}
\label{IM}
C:~V_z\longrightarrow V^{*b}_z~~~~, ~~v_i\otimes z^n\longrightarrow
v_{-i}^*\otimes z^n.
\end{eqnarray}

Throughout, we denote by $V(\lambda)$ a level-one irreducible
highest weight $\A$-module with 
highest weight $\lambda$. Consider the following intertwiners of 
$\A$-modules:
\begin{eqnarray*}
& &\Phi_{\lambda}^{\mu V}(z) :
 V(\lambda) \longrightarrow V(\mu)\otimes V_{z} ,\ \ \ \ 
\Phi_{\lambda}^{\mu V^{*}}(z) :
 V(\lambda) \longrightarrow V(\mu)\otimes V_{z}^{*b} ,\\
& &\Psi_{\lambda}^{V \mu}(z) :
 V(\lambda) \longrightarrow V_{z}\otimes V(\mu),\ \ \ \
\Psi_{\lambda}^{V^* \mu}(z) :
 V(\lambda) \longrightarrow V_{z}^{*b}\otimes V(\mu).
\end{eqnarray*}
They are intertwiners in the sense that for any $x\in \A$, 
\begin{eqnarray}
\label{EX}
\Theta(z)\cdot x=\Delta(x)\cdot \Theta(z),\ \ \ \Theta(z)=
\Phi(z),\Phi^{*}(z),\Psi(z),\Psi^{*}(z).
\end{eqnarray}
$\Phi(z)$ ($\Phi^{*}(z)$) is called type I (dual) vertex operator
and $\Psi(z)$
($\Psi^{*}(z)$) type II (dual) vertex operator.

We expand the vertex operators as 
\begin{eqnarray}
& &\Phi(z)=\sum_{j=1,0,-1}\Phi_j(z)\otimes v_j\  ,\ \ \ \
\Phi^{*}(z)=\sum_{j=1,0,-1}\Phi^{*}_j(z)\otimes v^{*}_j,\\
& &\Psi(z)=\sum_{j=1,0,-1}v_j\otimes\Psi_j(z)\  ,\ \ \ \
\Psi^{*}(z)=\sum_{j=1,0,-1}v^{*}_j\otimes\Psi^{*}_j(z).
\end{eqnarray}
Define the operators $\Phi_j(z),\Phi^{*}_j(z),\Psi_j(z)$ and 
$\Psi^{*}_j(z)$ $(j=1,0,-1)$ acting on the Fock space $F_1$ by  
\bea
\label{VOX1}
& &\Phi_{-1}(z)=
exp\{\sum_{n=1}^{\infty}\frac{[2n]_{q_1}}{n}q^{\frac{9}{2}n}
\omega_{-n}(-z)^n\}
exp\{\sum_{n=1}^{\infty}\frac{[2n]_{q_1}}{n}q^{-\frac{7}{2}n}
\omega_{n}(-z)^{-n}\}e^{Q}(zq^4)^{P+\frac{1}{2}},\\
& & \Phi_{0}(z)=\{\frac{q^{-\frac{1}{2}}}{[2]_{q_1}}\}^{1\over 2}
[\Phi_{-1}(z),f_1]_q~~,~~
\Phi_{1}(z)=\{\frac{q^{1\over 2}}{[2]_{q_1}}\}^{1\over 2}
[\Phi_{0}(z),f_1]~~,\\
& &\Psi_{1}(z)=
exp\{-\sum_{n=1}^{\infty}\frac{[2n]_{q_1}}{n}q^{\frac{1}{2}n}
\omega_{-n}z^n\}
exp\{-\sum_{n=1}^{\infty}\frac{[2n]_{q_1}}{n}q^{-\frac{3}{2}n}
\omega_{n}z^{-n}\}e^{-Q}(-zq)^{-P+\frac{1}{2}},\\
& & \Psi_{0}(z)=\{\frac{q^{-1\over 2}}{[2]_{q_1}}\}^{1\over 2}
[\Psi_{1}(z),e_1]_q~~,~~
\Psi_{1}(z)=\{\frac{q^{1\over 2}}{[2]_{q_1}}\}^{1\over 2}
[\Psi_{0}(z),e_1]~~,\\
\label{VOX2}
& &\Phi_i^*(z)=\Phi_{-i}(-zq^{-3})~~,\Psi_i^*(z)=\Psi_{-i}(-zq^{3}),
\eea
where $[a,b]_x=ab-xba$ and
the bosonic $q$-oscillators $\omega_m$ are defined by
\begin{eqnarray*}
\omega_m=-\frac{m}{\{[4m]_{q_1}-(-1)^m[2m]_{q_1}\}
[2m]_{q_1}}a_m~~,
\end{eqnarray*}
\noindent such that $[a_n,\omega_m]=\delta_{m+n,0}$.

According to \cite{Jing}, the operators $\Phi(z),\;\Phi^*(z),\;
\Psi(z)$ and $\Psi^*(z)$ defined in (\ref{VOX1}---\ref{VOX2}) are the only 
vertex operators of $\A$ which intertwine the level-one irreducible
highest weight
$\A$-modules :
\begin{eqnarray}
& &\Phi(z) :
 V(\Lambda_1) \longrightarrow V(\Lambda_1)\otimes V_{z} ,\ \ \ \ 
\Phi^{*}(z) :
 V(\Lambda_1) \longrightarrow V(\Lambda_1)\otimes V_{z}^{*b} ,\no\\
& &\Psi(z) :
 V(\Lambda_1) \longrightarrow V_{z}\otimes V(\Lambda_1),\ \ \ \
\Psi^{*}(z) :
 V(\Lambda_1) \longrightarrow V_{z}^{*b}\otimes V(\Lambda_1).
\end{eqnarray}

\section{Exchange relations of vertex operators}
In this section, we derive the exchange relations of the type I and type II 
vertex operators of $\A$.

\subsection{The R-matrix}
We introduce some abbreviations:
\begin{eqnarray*}
& &(z;p_1,p_2,...,p_m)=\prod_{\{l_1,l_2,...,l_m\}=0}^{\infty}
(1-zp_1^{l_1}p_2^{l_2}
\cdot\cdot\cdot p_m^{l_m})~~,\\
& &\Theta_p(z)=(z;p)(pz^{-1};p)(p;p)~~,~~\{z\}=(z;q^6,q^6)~~.
\end{eqnarray*}

Let $\R(z)\in End(V\otimes V)$
be the R-matrix of $\A$: 
\begin{eqnarray}
\label{R}
\R(z)(v_i\otimes v_j)=\sum_{k,l}\R^{i,j}_{kl}(z)v_k\otimes v_l, \ \ \ \
 \forall v_i, v_j, v_k, v_l\in V,
\end{eqnarray}
where the matrix elements are given by
\begin{eqnarray*}
& &\R^{1,1}_{1,1}(\frac{z_1}{z_2})=\R^{-1,-1}_{-1,-1}(\frac{z_1}{z_2})=1,~~
\R^{-1,0}_{0,-1}(\frac{z_1}{z_2})=\R^{0,1}_{1,0}(\frac{z_1}{z_2})
=\frac{(q-q^{-1})z_2}{z_1q-z_2q^{-1}},\\
& &\R^{-1,0}_{-1,0}(\frac{z_1}{z_2})=\R^{0,-1}_{0,-1}(\frac{z_1}{z_2})
=\R^{1,0}_{1,0}(\frac{z_1}{z_2})=\R^{0,1}_{0,1}(\frac{z_1}{z_2})=
\frac{z_1-z_2}{z_1q-z_2q^{-1}},\\
& &\R^{0,-1}_{-1,0}(\frac{z_1}{z_2})=\R^{1,0}_{0,1}(\frac{z_1}{z_2})
=\frac{(q-q^{-1})z_1}{z_1q-z_2q^{-1}},~~
\R^{-1,1}_{-1,1}(\frac{z_1}{z_2})=
\R^{1,-1}_{1,-1}(\frac{z_1}{z_2})
=\frac{(z_1-z_2)q^2(z_1q+z_2)}{(z_1q-z_2q^{-1})(z_1q^4+z_2q)},\\
& &\R^{1,-1}_{-1,1}(\frac{z_1}{z_2})
=\frac{(q^2-1)z_1\{z_2+z_2q^2+z_1q^3-z_1q^2\}}
{(z_1q-z_2q^{-1})(z_1q^4+z_2q)},~~
\R^{-1,1}_{1,-1}(\frac{z_1}{z_2})
=\frac{(q^2-1)z_2\{z_2-z_2q+z_1q^3+z_1q\}}
{(z_1q-z_2q^{-1})(z_1q^4+z_2q)},\\
& &\R^{1,-1}_{0,0}(\frac{z_1}{z_2})=
\R^{0,0}_{-1,1}(\frac{z_1}{z_2})
=\frac{(q-q^{-1})q^{\frac{7}{2}}z_1(z_1-z_2)}
{(z_1q-z_2q^{-1})(z_1q^4+z_2q)},\\
& &\R^{-1,1}_{0,0}(\frac{z_1}{z_2})=
\R^{0,0}_{1,-1}(\frac{z_1}{z_2})
=\frac{(q^{-1}-q)q^{\frac{3}{2}}z_2(z_1-z_2)}
{(z_1q-z_2q^{-1})(z_1q^4+z_2q)},\\
& &\R^{0,0}_{0,0}(\frac{z_1}{z_2})
=\frac{z_1-z_2}{z_1q-z_2q^{-1}}+
\frac{(q-q^{-1})(q+q^4)z_2z_1}
{(z_1q-z_2q^{-1})(z_1q^4+z_2q)},\\
& &R^{ij}_{kl}=0\ \ ,\ \ {\rm otherwise}.
\end{eqnarray*} 
Define the R-matrices $R^{I}(z)$, $R^{II}(z)$ and $R^{U}(z)$ :
\begin{eqnarray}
\label{SC}
R^{I}(z)=r(z)\R(z)~~,~~R^{II}(z)=\overline{r}(z)\R(z)~~,~~
R^{U}(z)=\rho(z)\R(z)~~,
\end{eqnarray}
\noindent where
\begin{eqnarray*}
& &r(z)=z^{-1}\frac{(-q^3z;q^6)(q^2z;q^6)(q^6z^{-1};q^6)(-q^5z^{-1};q^6)}
{(-q^3z^{-1};q^6)(q^2z^{-1};q^6)(q^6z;q^6)(-q^5z;q^6)}~~,\\
& &\overline{r}(z)=\frac{(q^4z;q^6)(-q^3z;q^6)(-qz^{-1};q^6)(q^6z^{-1};q^6)}
{(-q^3z^{-1};q^6)(q^4z^{-1};q^6)(q^6z;q^6)(-qz;q^6)}~~,\\
& &\rho(z)=\frac{r(z)}{\tau(zq^{-1})}~~,~~{\rm and }~~
\tau(z)=z^{-1}\frac{\Theta_{q^6}(qz)\Theta_{q^6}(-q^2z)}
{\Theta_{q^6}(q^5z)\Theta_{q^6}(-q^4z)}.
\end{eqnarray*}

Then the R-matrices satisfy  Yang-Baxter equation(YBE) on
$V\otimes V\otimes V$
\begin{eqnarray*}
R^{i}_{12}(z)R^{i}_{13}(zw)R^{i}_{23}(w)=R^{i}_{23}(w)R^{i}_{13}(zw)
R^{i}_{12}(z),~~{\rm for}~~ i=I, II, U ~~~~,
\end{eqnarray*}
and moreover enjoys: (i) the initial condition,
$R^{i}(1)=P$ for $i=I,II,U $,  with $P$ being
the permutation operator; (ii) the unitarity condition,
\begin{eqnarray*}
R^{i}_{12}(\frac{z}{w})R^{i}_{21}(\frac{w}{z})=1,~~{\rm for }~~
i=I,II,U
\end{eqnarray*}
\noindent where $R^{i}_{21}(z)=PR^{i}_{12}(z)P$; and (iii) the crossing
relations
\begin{eqnarray*}
& &(R^{i})^{k,l}_{m,n}(z)=(-q^{\frac{1}{2}})^{-\rho_l}(R^{(i)})^{-n,k}_{-l,m}
(-z^{-1}q^{-3})(-q^{\frac{1}{2}})^{\rho_n},~~{\rm for} ~~
i=I,II.
\end{eqnarray*}
\noindent Here and throughout, 
\begin{eqnarray}
\label{CR}
\rho_1=-1~~,~~\rho_0=0~~,~~\rho_{-1}=1.
\end{eqnarray}

\subsection{The  Faddeev-Zamolodchikov algebra}
Now, we are in the position to calculate the exchange relations of the type I 
and type II vertex operators of $\A$ given in (\ref{VOX1}---\ref{VOX2}).

Firstly we bosonize  the derivation operator $d$:
\begin{eqnarray}
\label{DER}
d=\sum_{n=1}^{\infty}na_{-n}\omega_n-\frac{P^2}{2}~~.
\end{eqnarray}
\noindent One can easily check that this $d$ operator obeys the following
commutation relations :
\begin{eqnarray}
\label{DER1}
q^dX^{\pm}(z)q^{-d}=X^{\pm}(zq^{-1})~~,~~
q^d\Phi_i(z)q^{-d}=\Phi_i(zq^{-1})~~,~~
q^d\Psi_i(z)q^{-d}=\Psi_i(zq^{-1})~~,
\end{eqnarray}
\noindent as required.

Define
\begin{eqnarray*}
\oint dzf(z)=Res(f)=f_{-1} \ \ ,\ \ {\rm for \ \ formal\ \  series \ \ 
function\ \ } f(z)=\sum_{n\in Z}f_nz^{n} .
\end{eqnarray*}
\noindent Then the Chevalley generators of $\A$ can be expressed 
by the integrals,
\begin{eqnarray*}
& &e_1=\oint z^{-1}dzX^+(z),~~ 
  f_1=\oint z^{-1}dz X^-(z), \\
& &e_0=K^{-2}\oint\oint dzw^{-1}dw~[X^-(w),X^-(z)]_{q},\\
& &f_0=\frac{1}{[4]_{q_1}}\oint\oint
  z^{-2}dzw^{-1}dw~[X^+(z),X^+(w)]_{q^{-1}}.
\end{eqnarray*}
From the normal order relations in appendix A, one can also obtain the
integral expression of the vertex 
operators defined in (\ref{VOX1}---\ref{VOX2})
\bea
\Phi_0(z)&=&\{\frac{q^{-\frac{1}{2}}}{[2]_{q_1}}\}^{\frac{1}{2}}
\oint \frac{dw}{w}~\frac{(q^{-1}-q)}
{zq^3(1+\frac{zq^5}{w})(1+\frac{w}{zq^3})}
:\Phi_{-1}(z)X^{-}(w):,\\
\Phi_{1}(z)&=&\frac{1-q^2}{[2]_{q_1}}\oint \frac{dw_1}{w_1}
\oint \frac{dw}{w}~
\frac{(1-\frac{w}{w_1})\{w_1^2(1-\frac{w}{w_1q^2})
(1+\frac{wq}{w_1})(1+\frac{zq^5}{w_1})-zwq^2
(1+\frac{w_1}{zq^3})(1-\frac{wq^2}{w_1})(1+\frac{w_1q}{w})\}}
{z^3q^{12}(1+\frac{zq^5}{w_1})(1+\frac{wq}{w_1})
(1+\frac{w_1}{zq^3})
(1+\frac{w_1q^5}{w})
(1+\frac{w}{zq^5})
(1+\frac{w}{zq^3})}\no\\
& &~~\times :\Phi_{-1}(z)X^-(w)X^-(w_1):,\\
\Psi_0(z)&=&\{\frac{q^{-\frac{1}{2}}}{[2]_{q_1}}\}^{\frac{1}{2}}
\oint \frac{dw}{w}~\frac{(1-q^2)}
{wq(1-\frac{w}{zq^2})(1-\frac{z}{w})}
:\Psi_{1}(z)X^{+}(w):,\\
\Psi_{-1}(z)&=&\frac{1-q^2}{[2]_{q_1}}\oint \frac{dw_1}{w_1}
\oint \frac{dw}{w}~
\frac{(1-\frac{w}{w_1})\{
(1-\frac{w_1}{wq^2})(1-\frac{z}{w_1})(1+\frac{w}{w_1q})w_1-
zq(1-\frac{w}{w_1q^2})(1-\frac{w_1}{zq^2})(1+\frac{w_1}{wq})\}}
{zwq(1-\frac{w}{zq^2})(1-\frac{z}{w})(1-\frac{w_1}{zq^2})(1+\frac{w_1}{wq})
(1-\frac{z}{w_1})(1+\frac{w}{w_1q})}\no\\
& &~~\times~:\Psi_{1}(z)X^+(w)X^+(w_1):.
\eea

By the technique proposed in \cite{AJ} , using  the above
integral expressions and the relations given in appendix A and appendix B,
and after tedious calcualtions , we can show that
the bosonic vertex operators defined in 
(\ref{VOX1}---\ref{VOX2}) satisfy the Faddeev-Zamolodchikov (ZF) 
algebra
\begin{eqnarray}
\label{ZF1}
& &\Phi_j(z_2)\Phi_i(z_1)=\sum_{kl}R^{I}(\frac{z_1}{z_2})^{k,l}_{i,j}
\Phi_k(z_1)\Phi_l(z_2) ,\\
& &\Psi^{*}_i(z_1)\Psi^{*}_j(z_2)=-\sum_{kl}R^{II}(\frac{z_1}{z_2})_{kl}^{ij}
\Psi^{*}_l(z_2)\Psi^{*}_k(z_1) ,\\
\label{ZF2}
& &\Psi^{*}_i(z_1)\Phi_j(z_2)=\tau(\frac{z_1}{z_2})
\Phi_j(z_2)\Psi^{*}_i(z_1)~~~.
\end{eqnarray}
And moreover, 
the bosonic vertex operators define
in (\ref{VOX1}---\ref{VOX2}) have the following invertibility relations,
\bea
\label{DEF1}
&&\Phi_i(z)\Phi^{*}_j(z)=g\delta_{ij}(-q^{\frac{1}{2}})^{\rho_i},
~~ g=\frac{(q^6;q^6)(-q^5;q^6)}{(-q^3;q^6)(q^2;q^6)},\\
&&\sum_{k}(-q^{\frac{1}{2}})^{-\rho_k}\Phi^{*}_k(z)\Phi_k(z)=g,
\eea
and 
\bea
&&\Psi_i(z_1)\Psi^{*}_j(z_2)=\frac{\overline{g}\delta_{ij}
(-q^{\frac{1}{2}})^{-\rho_i}}{1-\frac{z_2}{z_1}}~ + ~{\rm
regular~~ term}~~~~~ {\rm when}~ z_1\longrightarrow z_2,\\
&&\overline{g}=\frac{q(q^4;q^6)(-q^3;q^6)}{(-q;q^6)(q^6;q^6)}~~.\no
\eea
In the derivation of the above relations the following fact is helpful,
\begin{eqnarray*}
& &:\Phi_{-1}(z)\Phi_{-1}(-zq^{-3})X^{-}(zq^2)X^-(-zq^3)(zq^2)^{-1}
(-zq^3)^{-1}:=id,\\
& &:\Psi_{1}(z)\Psi_{1}(-zq^{3})X^{+}(zq^2)X^+(-zq^3)(zq^2)^{-1}
(-zq^3)^{-1}:=id.
\end{eqnarray*}

\subsection{Miki's construction of $\A$}

We generalize the Miki's construction to the twisted quantum affine algebra
$\A$  case.

Define
\begin{eqnarray*}
& &(L^+(z))^{j}_{i}=\Phi_i(zq^{\frac{1}{2}})
\Psi^{*}_j(zq^{-\frac{1}{2}}), \\ 
& &(L^-(z))^{j}_{i}=\Phi_i(zq^{-\frac{1}{2}})\Psi^{*}_j(zq^{\frac{1}{2}}).
\end{eqnarray*}
Then from the Faddeev-Zamolodchikov algebra in
(\ref{ZF1}---\ref{ZF2}) and the identity 
\begin{eqnarray*}
-\frac{r(z)}{\overline{r}(z)}=
\frac{\tau(zq^{-1})}{\tau(z^{-1}q^{-1})}
=\frac{\tau(zq)}{\tau(z^{-1}q)}~~,
\end{eqnarray*}
\noindent we can verify by straightforward computations that
the L-operators $L^{\pm}(z)$  give a
realization of the  Reshetikhin-Semenov-Tian-Shansky (RS) algebra\cite{RS}
 at level one in the twisted quantum affine algebra $\A$ 
\begin{eqnarray*}
& &R^U(\frac{z}{w})L^{\pm}_1(z)L^{\pm}_2(w)=
L^{\pm}_2(w)L^{\pm}_1(z)R^U(\frac{z}{w}),\\
& &R^U(\frac{z^+}{w^-})L^{+}_1(z)L^{-}_2(w)=
L^{-}_2(w)L^{+}_1(z)R^U(\frac{z^-}{w^+}) ,
\end{eqnarray*}
where $L^{\pm}_1(z)=L^{\pm}(z)\otimes 1$, 
$L^{\pm}_2(z)=1\otimes L^{\pm}(z)$ and $z^{\pm}=zq^{\pm\frac{1}{2}}$.
We remark that the Drinfeld bases (\ref{DRB1}---\ref{DRB2}) of the twisted
quantum affine algebra
 $\A$ can also be derived  by using the twisted RS algbera and the
 corresponding
 Gauss decomposition\cite{YZP}.

\section{The Izergin-Korepin  model }
In this section, we give a mathematical definition of the Izergin-Korepin 
 model on an infinite lattice.

\subsection{Space of states}
By means of the R-matrix (\ref{R}) of $\A$,
one can define the Izergin-Korepin model on the infinte lattice $\cdots
\otimes V\otimes V\otimes V\cdots$. Let $h$ be the operator on
$V\otimes V$ such that
\begin{eqnarray*}
& &PR(\frac{z_1}{z_2})=1+u h+\cdots,\ \ \ \ \ \ \ \ u
\longrightarrow 0 ,\\
& &\ \ \ \ \  P:{\rm the\ \ \ \ permutation \ \ operator },\ \
 e^{u}\equiv \frac{z_1}{z_2}.
\end{eqnarray*}
The Hamiltonian $H$ of the Izergin-Korepin  model is defined by\cite{MN}
\begin{eqnarray}
\label{HAM}
H=\sum_{l\in Z}h_{l+1,l}.
\end{eqnarray}
$H$ acts formally on the infinite tensor product,
\beq
\cdots V\otimes V\otimes V\cdots.\label{vvv}
\eeq
It can be easily checked that 
\begin{eqnarray*}
[U'_q(A^{(2)}_2), H]=0,
\end{eqnarray*}
where $U'_q(A^{(2)}_2)$ is the subalgebra of $\A$ with the
derivation operator $d$ being dropped. So $U'_q(A^{(2)}_2)$ plays
the role of infinite dimensional {\it non-abelian symmetries} of
the Izergin-Korepin model on the infinite lattice.
Since  the level-one vertex operators only exist between
$V(\Lambda_1)$ and itself ,  following \cite{DFJMN},
we can replace the infinite tensor product
(\ref{vvv}) by the level-0 $\A$-module,
\begin{eqnarray*}
{\bf H}={\rm Hom}(V(\Lambda_1),V(\Lambda_1))\cong V(\Lambda_1)\otimes
  V(\Lambda_1)^{*b},
\end{eqnarray*}
where $V(\Lambda_1)$ is level-one irreducible highest weight
$\A$-module and $V(\Lambda_1)^{*b}$ is the dual module of $V(\Lambda_1)$.
By  theorem 2, this homomorphism can be realized by applying the type
I vertex operators repeatedly.
So we shall make the (hypothetical) identification:
\begin{eqnarray*}
``{\bf \rm the \ \ space \ \ of \ \ physical \ \ states } " =
V (\Lambda_1)\otimes V(\Lambda_1)^{*b}.
\end{eqnarray*}
Namely, we take
\begin{eqnarray*}
{\bf H}\equiv End(V(\Lambda_1)) 
\end{eqnarray*}
as the space of states of the Izergin-Korepin 
model on the infinite lattice.
The left action of $\A$ on ${\bf H}$ is defined by
\begin{eqnarray*}
x.f=\sum x_{(1)}\circ f\circ b(x_{(2)}),~~~
 \forall x\in\A,~f\in {\bf H},
\end{eqnarray*}
where we have used notation $\Delta(x)=\sum x_{(1)}
\otimes x_{(2)}$. A linear operator of the form $O=A\otimes B
$ ($A\in End(V(\Lambda_1))$ and $B\in End(V(\Lambda_1)^{*b})$ )
operate on a state $f$ as :$f\longrightarrow A\circ f \circ B^{t}$.

Note that
${\bf H}$ has the unique canonical element which is refered as to
 the physical vacuum \cite{JM} and denoted  by
$|vac>$
\begin{eqnarray}
\label{VAC}
|vac>=\chi^{-\frac{1}{2}}(-q)^{-\hat{\rho}},
\end{eqnarray}
\noindent where $\chi$ coincide with the character of $\A$-module
$V(\Lambda_1)$
\begin{eqnarray*}
\chi=tr_{V(\Lambda_1)}(q^{-2\hat{\rho}})=q^{\frac{1}{4}}(1+q)
\prod_{n=1}^{\infty}(1+q^{6n-1})(1+q^{6n+1}).
\end{eqnarray*}
\noindent The proof of the above  character formula is given in appendix C.

\subsection{Local structure and local operators}
Following Jimbo and Miwa \cite{JM}, we use the type I 
vertex operators and their variants to incorporate the local structure
into the space of physical states ${\bf H}$, that is to formulate the action
of local operators of the Izergin-Korepin model on
the infinite tensor product (\ref{vvv}) in terms of their actions on
${\bf H}$. 

Using the isomorphisms (c.f. theorem 2)
\bea
\Phi(z)&:&~V(\Lambda_1)\longrightarrow V(\Lambda_1)\otimes V_z,\no\\
\Phi^{*,t}(-zq^3)&:&~V_z\otimes V(\Lambda_1)^{*b}\longrightarrow
  V(\Lambda_1)^{*b},
\eea
were $t$ is the transposition on the quantum space,
we have the following identification:
\begin{eqnarray*}
V(\Lambda_1)\otimes V(\Lambda_1)^{*b}\stackrel{\sim}{\rightarrow}
V(\Lambda_1)\otimes V_z\otimes V(\Lambda_1)^{*b}\stackrel{\sim}{\rightarrow}
V(\Lambda_1)\otimes V(\Lambda_1)^{*b}.
\end{eqnarray*}
The resulting isomorphism can be identified with
the  translation (or shift) operator defined by
\begin{eqnarray*}
T=g^{-1}\sum_i\Phi_i(1)\otimes \Phi_i^{*,t}(-q^3).
\end{eqnarray*}
Its inverse is given by
\begin{eqnarray*}
T^{-1}=g^{-1}\sum_i\Phi_i^{*}(1)\otimes \Phi_i^{t}(-q^{-3}).
\end{eqnarray*}

Thus we can define the local operators on $V$ as operators on
${\bf H}$ \cite{JM}. Let us label the tensor components from the middle as
$1,2,\cdots$ for the left half and as $0,-1,-2,\cdots$ for
the right half. The operators acting on the site 1  are defined by
\begin{eqnarray}
E_{i,j}\stackrel{def}{=}E^{(1)}_{i,j}=
g^{-1}(-q^{-\frac{1}{2}})^{-\rho_j}\Phi^{*}_i(1)\Phi_j(1)\otimes id.
\end{eqnarray}
In particular, we have spin operator $S_z$
\begin{eqnarray*}
S_z=g^{-1}\{
\Phi^*_1(1)\Phi_1(1)(-q^{-\frac{1}{2}})
+\Phi^*_0(1)\Phi_0(1)
-\Phi^*_{-1}(1)\Phi_{-1}(1)(-q^{-\frac{1}{2}})^{-1}\}
\otimes id.
\end{eqnarray*}

More generally we set
\begin{eqnarray}
\label{L}
E^{(n)}_{i,j}=T^{-(n-1)}E_{i,j}T^{n-1}\ \ \ \ (n\in Z).
\end{eqnarray}
Then, from the invertibility relations of the type I vertex operators
of $\A$, we can show that
the local operators $E^{(n)}_{ij}$ acting on
${\bf H}$  satisfy the following relations:
\beq
E^{(m)}_{i,j}E^{(n)}_{k,l}=\left\{
  \begin{array}{ll}
  \delta_{jk}E^{(n)}_{i,l} & {\rm if}\ \ \ m=n \\
  E^{(n)}_{k,l}E^{(m)}_{i,j}&{\rm if}\ \ \ m\ne n.
  \end{array}
 \right. 
\eeq
Moreover, we can define transfer matrix of the Izergin-Korepin model on the
infinite lattice  as 
\begin{eqnarray}
\label{TRM}
T(z)=g^{-1}\sum_{i}\Phi_i(z)\otimes \Phi^{*,t}_i(-zq^3)~~.
\end{eqnarray}
\noindent The commutativity
 of the transfer matrix ,[T(z),T(w)]=0, 
follows from the ZF algebraic relations (\ref{ZF1}).
The translation operator is $T=T(1)$, as expected.
Comparing with the definition of Hamiltonian of the Izergin-Korepin model
(\ref{HAM}), the action of such a  Hamiltonian on {\bf H} can be
given by
\begin{eqnarray}
\label{NH}
H=z\frac{d}{dz}\{ ln T(z)\}|_{z=1}
\end{eqnarray}

As is expected from the physical point of view, 
the vacuum vector
$|vac>$ is translationally invariant and singlet (i.e. belong
to the trivial representation of $\A$)
\begin{eqnarray}
T|vac>=|vac> ,\no\\
x.|vac>=\epsilon(x)|vac>.
\end{eqnarray}
This is proved as follows.
\begin{eqnarray*}
T(z)|vac>
&=&\chi^{-\frac{1}{2}}g^{-1}\sum_i\Phi_i(z)(-q)^{-\hat{\rho}}
\Phi^{*}_i(-zq^3)\\
&\stackrel{def}{=}&\chi^{-\frac{1}{2}}g^{-1}\sum_i\Phi_i(z)
(-q)^{-3d}(-q^{\frac{1}{2}})^
{-\frac{1}{2}h_1}
\Phi^{*}_i(-zq^3)\\
&=&\chi^{-\frac{1}{2}}g^{-1}\sum_i\Phi_i(z)\Phi_i^{*}(zq^6)
(-q^{\frac{1}{2}})^{\rho_i}(-q)^{-\hat{\rho}}\\
&\stackrel{def}{=}&\chi^{-\frac{1}{2}}g^{-1}\sum_i
\Phi^*_{-i}(-zq^3)\Phi_{-i}(-zq^3)
(-q^{\frac{1}{2}})^{-\rho_{-i}}(-q)^{-\hat{\rho}}\\
&=&\chi^{-\frac{1}{2}}(-q)^{-\hat{\rho}}=|vac>,
\end{eqnarray*}
\noindent where we have used the fact $\rho_i=-\rho_{-i}$ ,
(\ref{DEF}) and (\ref{DEF1}).Similary,
\begin{eqnarray*}
x.|vac>&=&\chi^{-\frac{1}{2}}\sum x_{(1)}(-q)^{\hat{\rho}}b(x_{(2)})\\
&=&\chi^{-\frac{1}{2}}\sum x_{(1)}S(x_{(2)})(-q)^{\hat{\rho}}\\
&=&\epsilon(x)|vac>
\end{eqnarray*}
\noindent In the third line we have used  the axioms of Hopf algebra
\cite{JM}
\begin{eqnarray*}
m\circ(id\otimes S)\circ \Delta=\epsilon.
\end{eqnarray*}
\noindent This completes the proof.

For any local operator $O$ on ${\bf H}$ defined in (\ref{L}),
its vacuum expectation value
is given by
\begin{eqnarray*}                                                           
<vac|O|vac>\stackrel{def}{=}<O>=
  \frac{tr_{V(\Lambda_1)}(q^{-2\hat{\rho}}O)}
  {tr_{V(\Lambda_1)}(q^{-2\hat{\rho}})}=
  \frac{tr_{V(\Lambda_1)}(q^{-6d-\frac{1}{2}h_1}O)}
  {tr_{V(\Lambda_1)}(q^{-6d-\frac{1}{2}h_1})},
\end{eqnarray*}
where the normalization $<vac|vac>=1$ has been chosen. We shall
denote the correlator $<vac|O|vac>$ by $<O>$.

By the proposition 4 and the definition of the local operators $
E^{(n)}_{i,j}$ (\ref{L}), we have
\begin{eqnarray*}
<E^{(n)}_{i,j}>=<vac|T^{-(n-1)}E_{i,j}T^{(n-1)}|vac>=
<vac|E_{i,j}|vac>=<E_{i,j}>
\end{eqnarray*}

\subsection{The n-particle states and form factors}
In order to construct the general eigenstates of the transfer matrix of
the Izergin-Korepin model (\ref{TRM}), we employ the type
II vertex operators.Define the n-particle states
\begin{eqnarray}
\label{nst}
|\xi_n,\xi_{n-1},\cdot\cdot\cdot,\xi_1>
_{i_n,i_{n-1},\cdot\cdot\cdot,i_1}
&=&\overline{g}^{-\frac{n}{2}}
\Psi^*_{i_n}(\xi_n) \Psi^*_{i_{n-1}}(\xi_{n-1})
\cdot\cdot\cdot\Psi^*_{i_1}(\xi_1)|vac>\no\\
&=&\overline{g}^{-\frac{n}{2}}
\chi^{-\frac{1}{2}}\Psi^*_{i_n}(\xi_n) \Psi^*_{i_{n-1}}(\xi_{n-1})
\cdot\cdot\cdot\Psi^*_{i_1}(\xi_1)(-q)^{-\hat{\rho}},
\end{eqnarray}
\noindent and its dual states
\begin{eqnarray}
\label{nst1}
{}_{i_n,i_{n-1},\cdot\cdot\cdot,i_1}
<\xi_n,\xi_{n-1},\cdot\cdot\cdot,\xi_1|
=\overline{g}^{-\frac{n}{2}}
\chi^{-\frac{1}{2}}(-q)^{-\hat{\rho}}
\Psi_{i_1}(\xi_1) \Psi_{i_2}(\xi_2)
\cdot\cdot\cdot\Psi_{i_n}(\xi_n).
\end{eqnarray}
Using the commutation relation of the type I vertex operators
,(\ref{ZF2}), one can verify  that
\begin{eqnarray*}
T(z)|\xi_n,\cdot\cdot\cdot,\xi_1>_{i_n,i_{n-1},\cdot\cdot\cdot,i_1}
&=&g^{-1}\overline{g}^{-\frac{n}{2}}
\chi^{-\frac{1}{2}}\sum_i \Phi_i(z)
\Psi^*_{i_n}(\xi_n) \Psi^*_{i_{n-1}}(\xi_{n-1})
\cdot\cdot\cdot\Psi^*_{i_1}(\xi_1)(-q)^{-\hat{\rho}}\Phi^*_i(-zq^3)\\
&=&g^{-1}\overline{g}^{-\frac{n}{2}}
\chi^{-\frac{1}{2}}\sum_i \Phi_i(z)
\Psi^*_{i_n}(\xi_n) \Psi^*_{i_{n-1}}(\xi_{n-1})
\cdot\cdot\cdot\Psi^*_{i_1}(\xi_1)\Phi^*_i(zq^6)(-q^{\frac{1}{2}})^{\rho_i}
(-q)^{-\hat{\rho}}\\
&=&g^{-1}\overline{g}^{-\frac{n}{2}}
\chi^{-\frac{1}{2}}\sum_i \prod_{j=1}^{n}\tau(\xi_j/z)\Phi_i(z)
\Phi^*_i(zq^6)(-q^{\frac{1}{2}})^{\rho_i}
\Psi^*_{i_n}(\xi_n) \Psi^*_{i_{n-1}}(\xi_{n-1})\\
& &~~~~~~~\times\cdot\cdot\cdot\Psi^*_{i_1}(\xi_1)
(-q)^{-\hat{\rho}}\\
&=&\prod_{j=1}^{n}\tau(\xi_j/z)
|\xi_n,\cdot\cdot\cdot,\xi_1>_{i_n,i_{n-1},\cdot\cdot\cdot,i_1}.
\end{eqnarray*}
\noindent Here we have used $\tau(zq^6)=\tau(z)$.Therefore the n-particle
states (\ref{nst}) are the eigenstates of transfer matrix $T(z)$. Likewise
for its dual states (\ref{nst1}). Hence, the Hamiltonian $H$ (\ref{NH}) on
the n-particle states are given by
\begin{eqnarray*}
H|\xi_n,\cdot\cdot\cdot,\xi_1>_{i_n,i_{n-1},\cdot\cdot\cdot,i_1}
=\sum_{j=1}^{n}\epsilon(\xi_j)
|\xi_n,\cdot\cdot\cdot,\xi_1>_{i_n,i_{n-1},\cdot\cdot\cdot,i_1},
\end{eqnarray*}
\noindent where
\begin{eqnarray*}
\epsilon(z)=-z\frac{d}{z}ln\tau(z).
\end{eqnarray*}
\noindent These results coincide with the energy of the elementary 
excitations derived from the Bethe ansatz method\cite{VR}.

In much the same way as  the correlation functions, the form
factors of a local operator $O$  of the form $O=A\otimes B$
can be given by
\begin{eqnarray*}
<vac|O|\xi_n,\cdot\cdot\cdot,\xi_1>_{i_n,i_{n-1},\cdot\cdot\cdot,i_1}
=\frac{\overline{g}^{-\frac{n}{2}}tr_{V(\Lambda_1)}((-q)^{-\hat{\rho}}A
\Psi^*_{i_n}(\xi_n) \Psi^*_{i_{n-1}}(\xi_{n-1})
\cdot\cdot\cdot\Psi^*_{i_1}(\xi_1)
(-q)^{-\hat{\rho}}B^{t})}
{tr_{V(\Lambda_1)}(q^{-2\hat{\rho}})}
\end{eqnarray*}
In particular , for the spin operator $S_z$, we have
\begin{eqnarray*}
& &<vac|S_z|\xi_n,\cdot\cdot\cdot,\xi_1>_{i_n,i_{n-1},\cdot\cdot\cdot,i_1}\\
& &~~~~=g^{-1}\frac{\overline{g}^{-\frac{n}{2}}tr_{V(\Lambda_1)}
(q^{-2\hat{\rho}}
\{\Phi^{*}_1(1)\Phi_1(1)(-q^{\frac{1}{2}})
+\Phi^{*}_0(1)\Phi_0(1)
-\Phi^{*}_{-1}(1)\Phi_{-1}(1)(-q^{\frac{1}{2}})^{-1}\}
\Psi^*_{i_n}(\xi_n) \Psi^*_{i_{n-1}}(\xi_{n-1})
\cdot\cdot\cdot\Psi^*_{i_1}(\xi_1))}
{tr_{V(\Lambda_1)}(q^{-2\hat{\rho}})}.
\end{eqnarray*}
\section{Correlation functions} 
The aim of  this section is to calculate $<E_{mn}>$.
The generalization to the calculation of the  multi-point functions 
 and the form factors is straightforward.

Set 
\begin{eqnarray*}
P^m_n(z_1,z_2|q)=
\frac{tr_{V(\Lambda_1)}(q^{-2\hat{\rho}}
\Phi^{*}_m(z_1)\Phi_n(z_2))}
{tr_{V(\Lambda_1)}(q^{-2\hat{\rho}})} ,
\end{eqnarray*}
then $<E_{mn}>=g^{-1}P^m_n(z,z|q)(-q^{\frac{1}{2}})^{-\rho_n}$.
Using the Clavelli-Shapiro technique\cite{CS} which we will present it in
appendix C, we get 
\begin{eqnarray}
\label{P}
& &P^m_n(z_1,z_2|q)\stackrel{def}{=}\delta_{mn}P_m(z_1,z_2|q)\no\\
& &~~~~=\frac{\delta_{mn}(C_{-1})^2(C^{(-)})^2}{\chi}G(-\frac{z_2q^3}{z_1})
\oint \frac{dw}{w}\oint\frac{dw_1}{w_1}
\theta(\frac{-z_1z_2q^5}{ww_1})I_m(z_1,z_2;w,w_1),
\end{eqnarray}
\noindent  where $\theta (z)$ is the elliptic function,
\begin{eqnarray*}
\theta(z)\stackrel{def}{=}
\sum_{n\in Z}q^{3(n+\frac{1}{2})^2}z^n~~~,
\end{eqnarray*}
\begin{eqnarray*}
I_1(z_1,z_2;w,w_1)&=&
(\frac{-z_1z_2q^5}{ww_1})^{\frac{1}{2}}W_1(-\frac{wq^3}{z_1})
W_1(-\frac{w_1q^3}{z_1})
\{
\frac{w}{z_2q^4}W_1(\frac{w}{z_2})W_1(\frac{w_1}{z_2})H(\frac{w_1}{w})
+\frac{z_2q^5}{w}W_2(\frac{z_2}{w_1})W_2(\frac{z_2}{w})H(\frac{w}{w_1})\\
& & ~~~~~~~
-qW_2(\frac{z_2}{w})W_1(\frac{w_1}{z_2})H(\frac{w_1}{w})
-W_2(\frac{z_2}{w_1})W_1(\frac{w}{z_2})H(\frac{w}{w_1})\}~~~,\\
I_0(z_1,z_2;w,w_1)&=&
(\frac{-z_1z_2q^5}{ww_1})^{\frac{1}{2}}W_1(-\frac{wq^3}{z_1})
W_2(\frac{z_2}{w})H(\frac{w_1}{w})
\{\frac{z_1q^2}{w}W_2(-\frac{z_1}{wq^3})W_1(\frac{w_1}{z_2})
-\frac{z_2q^5}{w_1}W_1(-\frac{w_1q^3}{z_1})W_2(\frac{z_2}{w_1})\\
& & ~~~~~~~
+W_1(-\frac{wq^3}{z_1})W_1(\frac{w_1}{z_2})
+q^2(\frac{-z_1z_2q^5}{ww_1})^{\frac{1}{2}}
W_2(-\frac{z_1}{wq^3})W_2(\frac{z_2}{w_1})\}~~~,\\
I_{-1}(z_1,z_2;w,w_1)&=&
-(\frac{-z_1z_2q^5}{ww_1})^{\frac{3}{2}}W_2(\frac{z_2}{w})
W_2(\frac{z_2}{w_1})
\{
qW_2(-\frac{z_1}{wq^3})W_1(-\frac{w_1q^3}{z_1})H(\frac{w_1}{w})
+W_2(-\frac{z_1}{w_1q^3})W_1(-\frac{wq^3}{z_1})H(\frac{w}{w_1})\\
& &~~~~+\frac{w}{z_1q}W_1(-\frac{wq^3}{z_1})W_1(-\frac{w_1q^3}{z_1})
H(\frac{w_1}{w})
+\frac{z_1q^2}{w}W_2(-\frac{z_1}{w_1q^3})W_2(-\frac{z_1}{wq^3})
H(\frac{w}{w_1})\}
\end{eqnarray*}
\noindent and the functions $G(z)$, $W_1(z)$, $W_2(z)$, $H(z)$  and
the constant $C_{-1}$, $C^{(-)}$ are given in appendix D (\ref{CT1}---
\ref{CT2}).

In particular,the correlation function of spin operator $S_z$
 is given by 
\begin{eqnarray}
<vac|S_z|vac>
=P_1(z,z|q)(-q^{-\frac{1}{2}})+P_0(z,z|q)
-P_{-1}(z,z|q)(-q^{-\frac{1}{2}})^{-1}
\end{eqnarray}
\noindent where $P_m(z_1,z_2|q)$ are defined in (\ref{P}).

We now derive the difference equations satisfied by these one-point
functions.
Noting (\ref{DER1}) and using the cyclicity of trace,
we get the difference equations of $P_m(z_1,z_2|q)$
\begin{eqnarray*}
P_m(z_1,z_2|q)
&=&\chi^{-1}tr_{V(\Lambda_1)}(q^{-6d-\frac{1}{2}h_1}
\Phi^*_m(z_1)\Phi_m(z_2))\\
&=&\chi^{-1}tr_{V(\Lambda_1)}(q^{-6d-\frac{1}{2}h_1}
\Phi_{-m}(-z_1q^{-3})\Phi_m(z_2))\\
&=&\chi^{-1}\sum_{n}R^I(z_2,-z_1q^{-3})^{n,-n}_{m,-m}
tr_{V(\Lambda_1)}(q^{-6d-\frac{1}{2}h_1}
\Phi_n(z_2)\Phi^*_n(z_1))\\
&=&\chi^{-1}
\sum_{n}q^{-\rho_n}R^I(z_2,-z_1q^{-3})^{n,-n}_{m,-m}
tr_{V(\Lambda_1)}
(\Phi_n(z_2q^6)q^{-6d-\frac{1}{2}h_1}
\Phi^*_n(z_1))\\
&=&\sum_{n}q^{-\rho_n}R^I(z_2,-z_1q^{-3})^{n,-n}_{m,-m}P_n(z_1,z_2q^6|q)~~.
\end{eqnarray*}

\section*{Acknowledgements}
W.L.Yang would like to thank Prof.Fan for  fruitful
discussions.Y.-Z.Zhang  thanks Australian Research
Council
 IREX programme for an Asia-Pacific Link Award and Institute of Modern
Physics of Northwest University for hospitality. The financial support
from the National Natural Science Foundation of China
and Australian Research
Council large, small, QEII fellowship
grants is also gratefully acknowledge.

\section*{Appendix A.}
In this apppendix, we give the normal order relations of fundmental bosonic
fields:
\begin{eqnarray*}
& &\Phi_{-1}(z)\Phi_{-1}(w)=zq^4g(w/z):\Phi_{-1}(z)\Phi_{-1}(w):~~,\\
& &\Phi_{-1}(z)X^{-}(w)=\frac{1}{zq^4(1+\frac{w}{zq^3})}
:\Phi_{-1}(z)X^-(w):~~,\\
& &X^-(w)\Phi_{-1}(z)=\frac{1}{w(1+\frac{zq^5}{w})}
:\Phi_{-1}(z)X^-(w):~~,\\
& &\Phi_{-1}(z)X^{+}(w)=(zq^4+w):\Phi_{-1}(z)X^+(w):X^+(w)\Phi_{-1}(z)~~,\\
& &\Psi_{1}(z)\Psi_{1}(w)=-zq\overline{g}(w/z):
\Psi_{1}(z)\Psi_{1}(w):~~,\\
& &\Psi_{1}(z)\Phi_{-1}(w)=-(zq)^{-1}h_1(w/z):
\Psi_{1}(z)\Phi_{-1}(w):~~,\\
& &\Phi_{-1}(z)\Psi_{1}(w)=(zq^4)^{-1}h_2(w/z):
\Phi_{-1}(z)\Psi_{1}(w):~~,\\
& &\Psi_{1}(z)X^{-}(w)=(w-zq):\Phi_{-1}(z)X^-(w):~~,~~
X^{-}(w) \Psi_{1}(z), \\
& &X^+(w)\Psi_{1}(z)=\frac{1}{w-z}
:\Psi_{1}(z)X^+(w):~~,\\
& &\Psi_{1}(z)X^+(w)=-\frac{1}{zq(1-\frac{w}{zq^2})}
:\Psi_{1}(z)X^+(w):~~,
\end{eqnarray*}
\noindent where
\begin{eqnarray*}
& &g(z)=\frac{(-q^3z;q^6)(q^2z;q^6)}{(q^6z;q^6)(-q^5z;q^6)}~~,~~
\overline{g}(z)=\frac{(-qz;q^6)(z;q^6)}{(q^4z;q^6)(-q^3z;q^6)},\\
& &h_1(z)=\frac{(-q^8z;q^6)(q^7z;q^6)}{(q^5z;q^6)(-q^4z;q^6)}~~,~~
h_2(z)=\frac{(-q^2z;q^6)(qz;q^6)}{(q^{-1}z;q^6)(-q^{-2}z;q^6)},
\end{eqnarray*}
\noindent and
\begin{eqnarray*}
& &X^+(z)X^+(w)=\frac{(z-w)(z-wq^{-2})}{z+wq^{-1}}:X^+(z)X^+(w):~~,\\
& &X^-(z)X^-(w)=\frac{(z-w)(z-wq^{2})}{z+wq}:X^-(z)X^-(w):~~,\\
& &X^+(z)X^-(w)=\frac{z+w}{(z-wq)(z-wq^{-1})}:X^+(z)X^+(w):~~,\\
& &X^-(w)X^+(z)=\frac{z+w}{(w-zq)(w-zq^{-1})}:X^+(z)X^+(w):~~.
\end{eqnarray*}

\section*{Appendix B. }
By means of the bosonic realiztion of $\A$, the integral expressions of
the vertex operators and the technique given in Ref.\cite{AJ}, one 
can check the following relations

\begin{itemize}
\item For the type I vertex operators
\begin{eqnarray*}
& &[\Phi_1(z),f_1]_{q^{-1}}=0,\ \ \Phi_1(z)=\frac{1}{\alpha}
[\Phi_0(z),f_1],\ \
\Phi_0(z)=\frac{\alpha}{[2]_{q_1}}[\Phi_{-1}(z),f_1]_q,\\
& &[\Phi_1(z),e_1]=\frac{[2]_{q_1}}{\alpha}t_1\Phi_0(z),
\ \ [\Phi_0(z),e_1]=\alpha t_1\Phi_{-1}(z),\ \
[\Phi_{-1}(z),e_1]=0,\\
& &\Phi_1(z)t_1=qt_1\Phi_1(z),\ \
\Phi_0(z)t_1=t_1\Phi_0(z),\ \
\Phi_{-1}(z)t_1=q^{-1}t_1\Phi_{-1}(z),\\
\end{eqnarray*}
where we take $\alpha=\{[2]_{q_1}q^{-\frac{1}{2}}\}^{\frac{1}{2}}$.

\item For the type II vertex operators
\begin{eqnarray*}
& &\Psi_{-1}(z)=\frac{1}{\alpha}[\Psi_0(z),e_1],
\ \ \Psi_0(z)=\frac{\alpha}{[2]_{q_1}}[\Psi_1(z),e_1]_q,\ \
[\Psi_{-1}(z),e_1]_{q^{-1}}=0,\\
& &[\Psi_1(z),f_1]=0,\ \ [\Psi_0(z),f_1]=\alpha t_1^{-1}\Psi_1(z),\ \
[\Psi_{-1}(z),f_1]=\frac{[2]_{q_1}}{\alpha}t_1^{-1}\Psi_0(z),\\   
& &\Psi_1(z)t_1=qt_1\Psi_1(z),\ \ \Psi_0(z)t_1=t_1\Psi_0(z),
\ \
\Psi_{-1}(z)t_1=q^{-1}t_1\Psi_{-1}(z).
\end{eqnarray*}
\end{itemize}

\section*{Appendix C.}
In computating  the correlation 
functions, one encounters the trace of the form
\begin{eqnarray*}
tr(x^{-d}e^{\sum_{m=1}^{\infty}A_{m}a_{-m}}
e^{\sum_{m=1}^{\infty}B_{m}a_{m})}f^{P})~~~~,   
\end{eqnarray*}
where $A_{m}$ ,$ B_m$ and $f$ are all some
coefficients. 
We can calculate the contributions 
from the oscillators modes and the zero modes {\it separately}.
The trace over the oscillator modes can be carried out as follows 
by using the Clavelli-Shapiro 
technique \cite{CS}.
Let us introduce the extra oscillators 
$a'_m$ which commutate with $a_m$. 
$a'_m$ satisfy the same commutation relations as those 
satisfied by $a_m$. 
Introduce the operators 
\begin{eqnarray*}
& &H_m=\frac{a_m\otimes 1}{1-x^m}+1\otimes a'_{-m}~~~,~~~m>0~~~,\\ 
& &H_m=a_m\otimes 1+\frac{1\otimes a'_{-m}}{x^m-1}~~~~,~~~m<0~~~, 
\end{eqnarray*}
\noindent which act on the space of the tensor Fock spaces of $\{a_m\}$ and
$\{a'_m\}$.
Then for any bosonic operator $O(a_m)$, one can show
\begin{eqnarray*}
tr(x^{-d}O(a_m))=\frac{
<0|O(H_m)|0>
}{\prod _{n=1}(1-x^n)} 
\end{eqnarray*}
providing that $d$ satisfies the derivation properties (\ref{DRB1}).
We write $<0|O(H_m)|0>\equiv <<O(a_m)>>$. Then by the
Wick theorem, one obtains
\begin{eqnarray*}
& &<<\Phi_{-1}(z)\Phi_{-1}(w)>>=C_{-1}C_{-1}G(w/z)~~,\\
& &<<\Phi_{-1}(z)X^-(w)>>=C_{-1}C^{(-)}W_1(w/z)~~,\\
& &<<X^-(z)\Phi_{-1}(w)>>=C_{-1}C^{(-)}W_2(w/z)~~,\\
& &<<X^-(z)X^-(w)>>=C^{(-)}C^{(-)}H(w/z)~~,
\end{eqnarray*}
\noindent where 
\begin{eqnarray}
\label{CT1}
& &C_{-1}=\frac{\{-q^9\}\{q^8\}}{\{q^{12}\}\{-q^{11}\}}~~~,~~~
C^{(-)}=\frac{(q^6;q^6)(q^8;q^6)}{(-q^7;q^6)}~~,\\
& &G(z)=\frac{ \{-q^3z\} \{q^2z\}\{-q^9z^{-1}\}\{q^8z^{-1}\}}
{ \{q^6z\} \{-q^5z\}\{q^{12}z^{-1}\}\{-q^{11}z^{-1}\}}~~~,\\
& &W_1(z)=\frac{(-q^{11}z^{-1};q^6)}{(-q^{-3}z;q^6)}~~~,~~
W_2(z)=\frac{(-q^{3}z^{-1};q^6)}{(-q^{5}z;q^6)}~~~,~~\\
\label{CT2}
& &H(z)=\frac{(q^2z;q^6)(z;q^6)(q^8z^{-1};q^6)(q^6z^{-1};q^6)}
{(-qz;q^6)(-q^7z^{-1};q^6)}~~~.
\end{eqnarray}

Now , we use the above technique to calculate the character of
$\A$-module $V(\Lambda_1)$ 
\begin{eqnarray*}
\chi&=&tr_{V(\Lambda_1)}(q^{-2\hat{\rho}})=
tr_{F_1}(q^{-6d}q^{-\frac{1}{2}h_1})
=tr_{F_1}(q^{-6d}q^{-P})\\
&=&\frac{1}{\prod_{n=1}^{\infty}(1-q^{6n})}
\sum_{n\in Z}q^{3(n+\frac{1}{2})^2}q^{-(n+\frac{1}{2})}\\
&=&q^{\frac{1}{4}}(1+q)\prod_{n=1}^{\infty}(1+q^{6n-1})(1+q^{6n+1})~~~.
\end{eqnarray*}
\noindent We have used the Jacobi tripe product identity\cite{Mu}
\begin{eqnarray*}
\sum_{n\in Z}q^{(n+v-\frac{1}{2})^2}t^n
=q^{\frac{1}{12}}q^{v^2-v+\frac{1}{6}}
\prod_{n=1}^{\infty}(1-q^{2n})(1+q^{2(n+v-1)}t) (1+q^{2(n-v)}t^{-1})
\end{eqnarray*}

\end{document}